\documentclass[12pt]{article}
\usepackage{graphicx}
\usepackage[T1]{fontenc}
\usepackage{lmodern}
\usepackage{mathtools}
\usepackage{microtype}
\usepackage{xcolor}
\usepackage{enumitem}
\usepackage{hyperref}
\usepackage{amsmath,amsthm,amssymb,amsfonts,euscript,enumerate}
\usepackage{color,wrapfig}
\usepackage{pxfonts}
\usepackage{verbatim}
\numberwithin{equation}{section}
\newtheorem{theorem}{Theorem}[section]
\newtheorem{lemma}[theorem]{Lemma}
\newtheorem{proposition}[theorem]{Proposition}

\theoremstyle{definition}

\theoremstyle{remark}

\newcommand{\R}{\mathbb R}

\newcommand{\dd}{\,\mathrm d}

\newcommand{\loc}{\mathrm{loc}}
\newcommand{\norm}[1]{\left\lVert #1\right\rVert}

\newcommand{\Cstar}{C_*}
\newcommand{\dist}{\operatorname{dist}}

\newcommand{\3}{\varepsilon}

\begin{document}
\title{Convergence of solutions of the fast diffusion equation to the singular eternal self-similar profiles}
\author{Shu-Yu Hsu\\
Department of Mathematics\\
National Chung Cheng University\\
168 University Road, Min-Hsiung\\
Chia-Yi 621, Taiwan, R.O.C.\\
e-mail: shuyu.sy@gmail.com}
\date{Nov 29, 2026}
\smallbreak \maketitle
\begin{abstract}
Let $n\ge 3$, $0<m<\frac{n-2}{n}$,  $\beta<0$, and $\alpha=\frac{2\beta}{1-m}$. For any $\lambda>0$, let  $U_{\lambda}(x,t)=e^{-\alpha t}f_{\lambda}(e^{-\beta t}x)$ be the singular eternal self-similar solutions of $u_t=\frac{n-1}{m}\Delta u^m$ in $(\mathbb{R}^n\setminus\{0\})\times\mathbb{R}$ whose profile satisfies 
\[
 \lim_{r\downarrow0}\frac{r^2f_\lambda(r)^{1-m}}{\log(1/r)}
 =\frac{2(n-1)(n-2-nm)}{|\beta|(1-m)},\qquad
 \lim_{r\to\infty}r^{(n-2)/m}f_\lambda(r)
 =\lambda^{2/(1-m)-(n-2)/m}.
\] 
We extend the convergence theory of Hui--Park~\cite{HP} to the full
stated parameter range and to initial data without radial symmetry
or a positive lower profile bound.
We prove existence and uniqueness, within the class
$0<u\le U_{\bar\lambda}$, for initial data satisfying
\[
 0\le u_0\le f_{\bar\lambda}\quad\hbox{a.e.},\qquad
 \int_{\mathbb R^n\setminus\{0\}}|u_0-f_{\lambda_0}|\,|x|^{-\mu}\dd x<\infty,
 \quad \mu=n-\frac{2}{1-m},
\]
for some constants $\bar\lambda,\lambda_0>0$. We also prove that the rescaled solution
$v(x,t)=e^{\alpha t}u(e^{\beta t}x,t)$  satisfies
\[
 \lim_{t\to\infty}\|v(\cdot,t)-f_{\lambda_0}\|_{L^1(|x|^{-\mu})}=0,
 \qquad
 \lim_{t\to\infty}\|v(\cdot,t)-f_{\lambda_0}\|_{C^2(K)}=0
\]
for every compact set $K\subset\mathbb R^n\setminus\{0\}$.
\end{abstract}

\vskip 0.2truein



\vskip 0.2truein
\setcounter{section}{0}

\section{Introduction}
\setcounter{equation}{0}
\setcounter{theorem}{0}

Let $n\ge 3$ and $0<m<\frac{n-2}{n}$. Recently there is a lot of study on the equation,
\begin{equation}\label{eq:fde}
 u_t=\frac{n-1}{m}\Delta u^m,\qquad u>0.
\end{equation}
by P. Daskalopoulos,  N. Sesum,  \cite{DS1}, \cite{DS2},  S.Y. Hsu \cite{Hs},  K.M. Hui \cite{H1}, \cite{H2}, etc. We refer the readers to the survey paper \cite{A} by D.G.~Aronson
and the books \cite{DK}, \cite{V}, by P. Daskalopoulos, C.E.~Kenig and J.L.~Vazquez on the recent results of \eqref{eq:fde}.

The equation \eqref{eq:fde} arises in many physical models and in the study of Yamabe flow. When $m>1$, \eqref{eq:fde} is called the porous medium equation which arises in the modelling of gases passing through porous media and oil passing through sand  \cite{A}.  When $m=1$, it is the heat equation with diffusivity $n-1$. When  $n\ge 3$ and $m=\frac{n-2}{n+2}$, \eqref{eq:fde} arises in the study of the Yamabe flow \cite{DS1}, \cite{DS2}. 
  
Let
\begin{equation}\label{eq:parameters}
 n\geq3,\qquad 0<m<\frac{n-2}{n},\qquad
 \beta<0,\qquad \alpha=\frac{2\beta}{1-m},
\end{equation}
 and 
$f_1$ be the unique radially symmetric solution of 
\begin{equation*}
\frac{n-1}{m}\Delta f^m+\alpha f+\beta x\cdot\nabla f=0, \quad f>0,\quad\mbox{ in }\mathbb{R}^n\setminus\{0\},
\end{equation*}
given by Theorem 1.1 of \cite{HP} which satisfies
\begin{equation*}
\lim_{r\to 0} \frac{r^{2}f(r)^{1-m}}{\log r^{-1}}=\frac{2(n-1)(n-2-nm)}{ (1-m)|\beta|}
\end{equation*}
and
\begin{equation*}
\lim_{r \to \infty} r^{\frac{n-2}{m}} f(r)=1.
\end{equation*}
Let
\begin{equation*}
 f_\lambda(x)=\lambda^{2/(1-m)}f_1(\lambda x).
\end{equation*}
Then \(f_\lambda\) satisfies
\begin{equation*}
 \frac{n-1}{m}\Delta f_\lambda^m
       +\alpha f_\lambda+\beta x\cdot\nabla f_\lambda=0,\quad\mbox{ in }\mathbb{R}^n\setminus\{0\},
\end{equation*}
and
\begin{equation}\label{eq:f-profile-asymptotic}
\lim_{r\downarrow0}\frac{r^2f_{\lambda}(r)^{1-m}}{\log r^{-1}}=\frac{2(n-1)(n-2-nm)}{|\beta|(1-m)},\qquad \lim_{r\to\infty}r^{\frac{n-2}{m}}f_{\lambda}(r)=\lambda^{\frac{2}{1-m}-\frac{n-2}{m}}.
\end{equation}
Let
\begin{equation}\label{eq:self-similar-formula}
U_{\lambda}(x,t)=e^{-\alpha t}f_{\lambda}(e^{-\beta t}x).
\end{equation}
Then $U_{\lambda}$ is the eternal self-similar solution of \eqref{eq:fde} in $(\mathbb{R}^n\setminus\{0\})\times\mathbb{R}$. In \cite{HP} K.M.~Hui and J.~Park (Theorem 1.4 of \cite{HP}) proved that for any $0\le u_0\in L^1_{loc}(\mathbb{R}^n\setminus\{0\})$ satisfying
\begin{equation}\label{eq:sandwich}
 f_{\lambda_1}(x)\leq u_0(x)\leq f_{\lambda_2}(x)
 \quad\hbox{for a.e. }x\in{\mathbb{R}^n\setminus\{0\}},\qquad \lambda_1>\lambda_2>0,
\end{equation}
there exists a unique solution $u$ of
\begin{equation}\label{eq:pde}
\left\{\begin{aligned}
&u_t=\frac{n-1}{m}\Delta u^m\quad\hbox{in }\mathbb{R}^n\setminus\{0\}\times(0,\infty)\\
&u(\cdot,0)=u_0\qquad\quad\hbox{in }\mathbb{R}^n\setminus\{0\}.\end{aligned}\right.
\end{equation}
which satisfies 
\begin{equation}\label{eq:upperAB}
u_t \le \frac{u}{(1-m) t} \quad \text{ in } (\mathbb{R}^n \setminus \{0\}) \times (0, \infty)
\end{equation}
and
\begin{equation}\label{eq:physicaltrap}
 U_{\lambda_1}(x,t)\leq u(x,t)\leq U_{\lambda_2}(x,t)\quad\forall x\in\mathbb{R}^n\setminus\{0\}, t>0.
\end{equation}
The rescaled function of $u$ is defined as
\begin{equation}\label{eq:rescaling}
 v(x,t):=e^{\alpha t}u(e^{\beta t}x,t)\quad\forall x\in\mathbb{R}^n\setminus\{0\}, t>0.
\end{equation}
It satisfies
\begin{equation}\label{eq:rescaledpde}
 v_t=\frac{n-1}{m}\Delta v^m+\alpha v+\beta x\cdot\nabla v\quad \text{ in } (\mathbb{R}^n \setminus \{0\}) \times (0, \infty)
\end{equation}
Note that \eqref{eq:physicaltrap} holds if and only if 
\begin{equation}\label{eq:rescaed_function_trap}
\qquad f_{\lambda_1}(x)\leq v(x,t)\leq f_{\lambda_2}(x)\quad\forall x\in\mathbb{R}^n\setminus\{0\}, t>0.
\end{equation}
In the paper \cite{HP} K.M.~Hui and J.~Park  proved the asymptotic large time  behavior of such singular solution $u$ when  $n=3,4$ and $\frac{n-2}{n+2}\le m<\frac{n-2}{n}$ holds. They also  proved the asymptotic large time  behaviour for the case $3\le n<8$, $1-\sqrt{2/n}\le m<\min\left(\frac{2(n-2)}{3n},\frac{n-2}{n+2}\right)$,  and $u(x,t)$ is radially symmetric in $x\in\mathbb{R}^n\setminus\{0\}$ for any $t>0$ under appropriate conditions on the initial value $u_0$.

In this paper we will extend the results of \cite{HP} and prove the existence, uniqueness and asymptotic behavior of solutions of \eqref{eq:pde}  with initial value $u_0$ satisfying 
\begin{equation}\label{eq:upperinitial}
 0\leq u_0(x)\leq f_{\bar\lambda}(x)
 \quad\hbox{for a.e. }x\in{\mathbb{R}^n\setminus\{0\}}
\end{equation}
and 
\begin{equation}\label{eq:mainhyp}
 E_0:=\int_{\mathbb{R}^n\setminus\{0\}} |u_0(x)-f_{\lambda_0}(x)|
                    |x|^{-\mu}\dd x<\infty,
 \qquad \mu=n-\frac{2}{1-m}.
\end{equation}
for some constants $\lambda_0\ge\bar\lambda>0$. For any positive weight function \(w\) on $\mathbb{R}^n\setminus\{0\}$, we let
\begin{equation*}
\norm{h}_{L^1(w)}=\int_{\mathbb{R}^n\setminus\{0\}} |h(x)|w(x)\dd x.
\end{equation*}
Let
\begin{equation*}
 a=\frac{2}{1-m},\qquad b=\frac{n-2}{m},\qquad
 \mu=n-a,\qquad \gamma_1=b-a=\frac{\mu}{m},\qquad
 \delta=a-2.
\end{equation*}
Then
\begin{equation*}
 2<a<n<b,\qquad 0<\mu<n-2,\qquad \delta>0,\qquad
 \mu+\delta=n-2,\qquad am=a-2.
\end{equation*}
In this paper we will prove the following results. The first theorem is an analogue of the convergence results in
\cite[Theorems~1.6 and~1.7]{HP} for the full range \eqref{eq:parameters}.
Its weighted conclusion uses $|x|^{-\mu}$, rather than the particular
weights in those theorems. The second theorem removes the lower profile
bound. Neither theorem assumes radial symmetry of the initial data.

\begin{theorem}[Nonradial convergence in the critical weighted norm]\label{thm:main}
Suppose \eqref{eq:parameters} holds and $u_0$ satisfies \eqref{eq:sandwich} and \eqref{eq:mainhyp} for some constants  \(\lambda_0>0\) and $\lambda_1>\lambda_2>0$.
Then \(\lambda_0\) is uniquely determined by this condition and necessarily
\begin{equation}\label{eq:lambda-ordering}
\lambda_2\leq\lambda_0\leq\lambda_1.
\end{equation}
Let $u$ be the unique solution of \eqref{eq:pde} given by Theorem 1.4 of \cite{HP} which satisfies \eqref{eq:upperAB} and \eqref{eq:physicaltrap}.
Let $v$ be given by \eqref{eq:rescaling}. Then
\begin{equation}\label{eq:mainL1}
\lim_{t\to\infty}\int_{\mathbb{R}^n\setminus\{0\}}
       |v(x,t)-f_{\lambda_0}(x)|\,|x|^{-\mu}\dd x=0
\end{equation}
and
\begin{equation}\label{eq:mainC2}
 \lim_{t\to\infty}\norm{v(\cdot,t)-f_{\lambda_0}}_{C^2(K)}=0
 \quad\hbox{for any compact set } K\subset\mathbb{R}^n\setminus\{0\}. 
\end{equation}
\end{theorem}

\begin{theorem}[removal of lower profile bound condition]
\label{thm:upperonly}
Suppose \eqref{eq:parameters} holds and $u_0$ satisfies \eqref{eq:upperinitial} and 
\eqref{eq:mainhyp} for some constants $\lambda_0>0$, $\bar\lambda>0$. The parameter \(\lambda_0\) in \eqref{eq:mainhyp} is unique and
necessarily \(\lambda_0\geq\bar\lambda\). 
Then there exists a positive
classical solution $u$ of \eqref{eq:pde} satisfying \eqref{eq:upperAB} and
\begin{equation}\label{eq:upperclass}
 0<u(x,t)\leq U_{\bar\lambda}(x,t)
 \quad\forall x\in{\mathbb{R}^n\setminus\{0\}},\ t>0.
\end{equation}
This solution is unique among all positive classical solutions with
the same $L^1_{\loc}$ initial trace satisfying \eqref{eq:upperclass}.
Let $v$ be given by \eqref{eq:rescaling}. Then \(v\) satisfies both
\eqref{eq:mainL1} and \eqref{eq:mainC2}.
Furthermore, for any compact set \(K\subset{\mathbb{R}^n\setminus\{0\}}\) and constant \(h>0\),
there exists a constant \(c_{K,h}>0\) such that
\begin{equation}\label{eq:upperlocalpositive}
 c_{K,h}\leq v(x,t)\leq f_{\bar\lambda}(x)
 \quad\forall x\in K,\ t\geq h.
\end{equation}
\end{theorem}

\subsection{Definitions}

For any $0\leq u_0\in L^1_{\loc}(\R^n\setminus\{0\}),$ we say that $u$ is a solution of \eqref{eq:pde}
if $u>0$ in $(\R^n\setminus\{0\})\times(0, \infty)$ satisfies \eqref{eq:fde} in $ \left(\R^n\setminus\{0\}\right)\times(0, \infty)$ in the classical sense and
\begin{equation*}
\| u(\cdot, t)-u_0\|_{L^1(K)}\to 0\quad\mbox{ as }t\to 0
\end{equation*}
for any compact set $K\subset \R^n\setminus\{0\}$. For any solution $u$ of \eqref{eq:fde} in $ \left(\R^n\setminus\{0\}\right)\times(0, \infty)$ we will also let $v$ be given by \eqref{eq:rescaling} for the rest of the paper.

 For any $d\in\mathbb{Z}^+$, $x_0\in\R^d$ and $R>0$, we let $B_R^d(x_0)=\{x\in\mathbb{R}^d: |x-x_0|<R\}$ and $B_R^d=B_R^d(0)$. When the dimension is $n$, we also write $B_R(x_0)=B_R^n(x_0)$. For any $R>1$, let $A_R=\{x\in\R^n:1/R<|x|<R\}$. We also let $\omega_d$ denote the surface area of
 $\partial B_1^d$. Then $|B_R^d|:=vol(B_R^d)=\omega_dR^d/d$. 
 
 The notation $K\Subset D$ means that $K$ is compact and contained in
the open set $D$. For a smooth function on a neighborhood of $K$,
$\|\cdot\|_{C^2(K)}$ is the maximum of the suprema on $K$ of its
spatial derivatives of order at most two.

\section{Convergence result for solutions bounded between two eternal solutions}
\setcounter{equation}{0}
\setcounter{theorem}{0}

In this section we will prove Theorem \ref{thm:main}. We first recall some results of \cite{HP}. Let
\begin{equation*}
 \Cstar=\frac{2(n-1)(n-2-nm)}{(1-m)|\beta|}>0,\qquad
 q=\frac{1}{1-m}.
\end{equation*}
Then by  \eqref{eq:f-profile-asymptotic} and Remark 3 of \cite{HP},  
for any \(\lambda>0\), 
\begin{equation}\label{eq:originasymp}
 f_\lambda(r)\sim r^{-a}
\bigl(\Cstar\log(1/r)\bigr)^q\quad\mbox{ as }r\downarrow0,
\end{equation}
\begin{equation}\label{eq:infinityasymp}
 f_\lambda(r)\sim\lambda^{-\gamma_1}r^{-b}\quad\mbox{ as }r\to\infty
\end{equation}
and
\begin{equation}\label{eq:ordering}
 \frac{\partial f_\lambda}{\partial\lambda}(r)<0\quad\forall r>0,\lambda>0.
\end{equation}
In particular, each \(f_\lambda\) is locally integrable across the origin,
because \(a<n\), and decays with exponent \(b>n\) at infinity.
The scaling formula also implies continuity of
\(\lambda\mapsto f_\lambda\) in \(C^2(K)\) for \(K\Subset{\mathbb{R}^n\setminus\{0\}}\).

The higher-order expansion \cite[Theorem~1.2, equation~(34)]{HP}
can be written, for each fixed $\lambda>0$, as
\[
 r^2 f_\lambda(r)^{1-m}
 =\Cstar\left(L+\kappa_0\log L+B_0-\log\lambda+o(1)\right),
 \qquad L=\log(1/r)\longrightarrow\infty,
\]
where
$\kappa_0=(n-2-(n+2)m)/(2(n-2-nm))$ and $B_0$ is independent of
$\lambda$ (the parameters $n,m,\beta$ are fixed).
Subtracting the expansions for two fixed parameters gives
\begin{equation}\label{eq:pressuredifference}
 \lim_{r\downarrow0}
 r^2\bigl(f_\lambda(r)^{1-m}
          -f_\nu(r)^{1-m}\bigr)
       =\Cstar\log\frac{\nu}{\lambda}
 \quad\forall\lambda,\nu>0.
\end{equation}

\begin{lemma}[Different profiles have infinite critical distance]
\label{lem:profiledistance}
For any $\lambda>0$, $\nu>0$ and $\lambda\ne\nu$,  we have
\begin{equation}\label{eq:profilegap}
 |f_\lambda(r)-f_\nu(r)|
 \sim \frac{\Cstar^q}{1-m}
       \left|\log\frac{\nu}{\lambda}\right|
       r^{-a}\bigl(\log(1/r)\bigr)^{m/(1-m)}
 \quad\mbox{ as }r\downarrow0.
\end{equation}
Consequently,
\begin{equation}\label{eq:infinitedistance}
 \int_{0<|x|<1}|f_\lambda(x)-f_\nu(x)|\,|x|^{-\mu}\dd x=\infty\quad\forall\lambda,\nu>0, \lambda\ne\nu.
\end{equation}
\end{lemma}

\begin{proof}
Let \(P_\lambda(r)=r^2f_\lambda(r)^{1-m}\). Then
\(f_\lambda(r)=r^{-a}P_\lambda(r)^q\). By \eqref{eq:originasymp},
\begin{equation}\label{eq:p-lambda-approx-value}
P_\lambda(r)\sim\Cstar\log(1/r)\quad\mbox{ as }r\to 0
\end{equation}
and by \eqref{eq:pressuredifference},
\begin{equation}\label{eq:p-lambda-difference}
P_\lambda(r)-P_\nu(r)\to\Cstar\log(\nu/\lambda)\quad\mbox{ as }r\to 0.
\end{equation}
By the mean value theorem for any $r>0$ there exists a constant $\xi(r)$ between $P_\lambda(r)$ and $P_\nu(r)$ such that
\begin{equation}\label{eq:p-lambda-q-difference}
|P_\lambda(r)^q-P_\nu(r)^q|=q\xi(r)^{q-1}|P_\lambda(r)-P_\nu(r)|\quad\forall r>0.
\end{equation}
By \eqref{eq:p-lambda-approx-value}, \eqref{eq:p-lambda-difference} and \eqref{eq:p-lambda-q-difference},
\begin{equation*}
 r^{\frac{2}{1-m}}|f_\lambda(r)-f_\nu(r)|
 \sim \frac{\Cstar^q}{1-m}
       \left|\log\frac{\nu}{\lambda}\right|
      \bigl(\log(1/r)\bigr)^{m/(1-m)}
 \quad\mbox{ as }r\downarrow0
\end{equation*}
and \eqref{eq:profilegap} follows. By \eqref{eq:profilegap} there exists a constant $0<\3<1$ such that
\begin{align*}
 |f_\lambda(r)-f_\nu(r)|\ge &C r^{-a}\bigl(\log(1/r)\bigr)^{m/(1-m)}
 \quad\forall 0<r<\3\\
 \Rightarrow\quad\int_{0<|x|<\3}|f_\lambda(x)-f_\nu(x)|\,|x|^{-\mu}\dd x\ge &C\int_0^{\3}r^{n-1}r^{-\mu}r^{-a}\bigl(\log(1/r)\bigr)^{m/(1-m)}\,dr\\
 =&C\int_0^{\3}r^{-1}\bigl(\log(1/r)\bigr)^{m/(1-m)}\,dr\\
 =&\infty
\end{align*}
and \eqref{eq:infinitedistance} follows.
\end{proof}

\begin{lemma}\label{lem:parameter-unique-and-range}
Suppose the hypotheses of Theorem \ref{thm:main} holds. Then the parameter \(\lambda_0\) in \eqref{eq:mainhyp} is unique and \(\lambda_2\leq\lambda_0\leq\lambda_1\). 
\end{lemma}
\begin{proof}
Suppose \eqref{eq:mainhyp} holds for two different parameters $\lambda_0\ne\bar\lambda_0$. Then 
\begin{align*}
\int_{\mathbb{R}^n\setminus\{0\}} |f_{\bar\lambda_0}(x)-f_{\lambda_0}(x)||x|^{-\mu}\dd x
\le&\int_{\mathbb{R}^n\setminus\{0\}} |u_0(x)-f_{\bar\lambda_0}(x)|
|x|^{-\mu}\dd x+\int_{\mathbb{R}^n\setminus\{0\}} |u_0(x)-f_{\lambda_0}(x)|
|x|^{-\mu}\dd x\\
<&\infty
\end{align*}
which contradicts Lemma \ref{lem:profiledistance}. Hence  the parameter $\lambda_0$ is 
unique.
 
Suppose \(\lambda_0>\lambda_1\). Then by \eqref{eq:sandwich} and Lemma \ref{lem:profiledistance},
\begin{align*}
&u_0-f_{\lambda_0}\geq f_{\lambda_1}-f_{\lambda_0}>0\quad\hbox{for a.e. }x\in\mathbb{R}^n\setminus\{0\}\\
\Rightarrow\quad &\int_{\mathbb{R}^n\setminus\{0\}} |u_0(x)-f_{\lambda_0}(x)|
|x|^{-\mu}\dd x\ge
\int_{\mathbb{R}^n\setminus\{0\}} |f_{\lambda_1}(x)-f_{\lambda_0}(x)|
|x|^{-\mu}\dd x=\infty
\end{align*}
which contradicts \eqref{eq:mainhyp}. Hence \(\lambda_0\le\lambda_1\). If \(\lambda_0<\lambda_2\),
Then by \eqref{eq:sandwich} and Lemma \ref{lem:profiledistance},
\begin{align*}
&f_{\lambda_0}-u_0\geq f_{\lambda_0}-f_{\lambda_2}>0\quad\hbox{for a.e. }x\in\mathbb{R}^n\setminus\{0\}\\
\Rightarrow\quad &\int_{\mathbb{R}^n\setminus\{0\}} |u_0(x)-f_{\lambda_0}(x)|
|x|^{-\mu}\dd x\ge
\int_{\mathbb{R}^n\setminus\{0\}} |f_{\lambda_0}(x)-f_{\lambda_2}(x)|
|x|^{-\mu}\dd x=\infty
\end{align*}
which contradicts \eqref{eq:mainhyp}. Hence \(\lambda_0\ge\lambda_2\). Thus \(\lambda_2\leq\lambda_0\leq\lambda_1\).
\end{proof}

For any \(s\geq0\), we let
\begin{equation*}
 W_s(r)=r^{-\mu}(1+sr)^{-\delta}
\end{equation*}
and
\begin{equation*}
 \Psi_{s,\theta}(r)=r^\theta W_s(r),
 \qquad\forall 0<\theta<\mu.
\end{equation*}
In particular \(W_0=r^{-\mu}\). When \(s>0\), 
\begin{equation*}
W_s(r)\sim r^{-\mu}\quad\mbox{ as }r\to 0
\end{equation*}
and
\begin{equation*}
W_s(r)\sim s^{-\delta}r^{2-n}\quad\mbox{ as }r\to \infty.
\end{equation*}

\begin{lemma}\label{lem:superharmonic}
For any \(s\geq0\) and \(0<\theta<\mu\),
\begin{equation}\label{eq:superharmonic}
 \Delta\Psi_{s,\theta}<0,\qquad
 |\nabla\Psi_{s,\theta}|
       \leq \frac{n-2}{r}\Psi_{s,\theta}
 \quad\hbox{on }{\mathbb{R}^n\setminus\{0\}}
\end{equation}
and
\begin{equation}\label{eq:Ws-superharmonic}
\Delta W_s\le 0\quad\hbox{on }{\mathbb{R}^n\setminus\{0\}}.
\end{equation}
\end{lemma}

\begin{proof}
Let 
\[
 k(r)=\mu-\theta+\frac{\delta sr}{1+sr}.
\]
Then \begin{equation}\label{eq:psi-derivative-eqn}
\Psi_{s,\theta}'/\Psi_{s,\theta}=-k(r)/r\quad\forall r>0,
\end{equation} 
and direct differentiation gives
\begin{align}
 \frac{\Psi_{s,\theta}''}{\Psi_{s,\theta}}
       &=\frac{k(r)^2+k(r)-rk'(r)}{r^2}\quad\forall r>0,\notag\\
 \frac{\Delta\Psi_{s,\theta}}{\Psi_{s,\theta}}
       &=\frac{k(r)(k(r)-n+2)-rk'(r)}{r^2}\quad\forall r>0,
 & rk'(r)&=\frac{\delta sr}{(1+sr)^2}\geq0\quad\forall r>0. \label{eq:laplacianweight}
\end{align}
Since \(\mu+\delta=n-2\),
\[
 0<\mu-\theta\leq k(r)\leq\mu+\delta-\theta<n-2.
\]
Since the first term in the numerator of
\eqref{eq:laplacianweight} is strictly negative and the second is
nonpositive, the first inequality of \eqref{eq:superharmonic} follows.
Since \(k<n-2\), by \eqref{eq:psi-derivative-eqn} the second inequality of \eqref{eq:superharmonic} follows. The same formula with \(\theta=0\) proves \eqref{eq:Ws-superharmonic}.
\end{proof}

\begin{proposition}\label{prop:contraction}
Let \(u_{0,1},u_{0,2}\) satisfy  
\begin{equation}\label{eq:initial-value-upper-bd}
0\leq u_{0,i}\leq f_{\bar\lambda}\quad\mbox{ in }\mathbb{R}^n\setminus\{0\}\quad\forall i=1,2
\end{equation}
for some constant $\bar\lambda>0$ and let \(u_1,u_2\) be the solutions of \eqref{eq:pde} with initial values \(u_0=u_{0,1},u_{0,2}\) respectively which satisfy 
\begin{equation}\label{eq:physicaltrap3}
0<u_i(x,t)\leq U_{\bar\lambda}(x,t)\quad\forall x\in\mathbb{R}^n\setminus\{0\}, t>0,i=1,2.
\end{equation}
For \(s\geq0\), if
\begin{equation}\label{eq:initial-data-weighted-difference}
\int_{\mathbb{R}^n\setminus\{0\}} |u_{0,1}-u_{0,2}|W_s\dd x<\infty,
\end{equation}
then
\begin{equation}\label{eq:contraction}
 \int_{\mathbb{R}^n\setminus\{0\}} |u_1(x,t)-u_2(x,t)|W_s(|x|)\dd x
 \leq \int_{\mathbb{R}^n\setminus\{0\}} |u_{0,1}-u_{0,2}|W_s\dd x
 \qquad\forall t>0.
\end{equation}
Similarly for \(s\geq0\), if
\begin{equation}\label{eq:initial-data-order}
\int_{\mathbb{R}^n\setminus\{0\}} (u_{0,1}-u_{0,2})_+W_s\dd x<\infty,
\end{equation}
then
\begin{equation}\label{eq:contraction-order}
 \int_{\mathbb{R}^n\setminus\{0\}} (u_1(x,t)-u_2(x,t))_+W_s(|x|)\dd x
 \leq \int_{\mathbb{R}^n\setminus\{0\}} (u_{0,1}-u_{0,2})_+W_s\dd x
 \qquad\forall t>0.
\end{equation}
\end{proposition}
\begin{proof}
We will use a modification of the proof of Theorem 1.3 of \cite{HP} and Lemma 4.1 of \cite{DS1} to prove the lemma.
By the Kato inequality (\cite{DK}, \cite{K}),
\begin{equation}\label{eq:kato}
\frac{\partial}{\partial t}|u_1-u_2|\leq \frac{n-1}{m}\Delta |u_1^m-u_2^m|
 \quad\hbox{in }\mathcal D'({\mathbb{R}^n\setminus\{0\}}\times(0,\infty)).
\end{equation}
We now fix \(T>0\). Since by \eqref{eq:physicaltrap3},
\[
 |u_1-u_2|\leq U_{\bar\lambda},\qquad
 |u_1^m-u_2^m|\leq U_{\bar\lambda}^m,
\]
by \eqref{eq:self-similar-formula},
\eqref{eq:originasymp}, \eqref{eq:infinityasymp} and \eqref{eq:physicaltrap3}, there exist constants
\(r_T\in(0,1/2)\), \(R_T>2\), and \(C_T>0\) such that, for
\(0<t\leq T\),
\begin{align}
 |u_1-u_2|(x,t)&\leq C_T r^{-a}(\log(1/r))^q,
 & |u_1^m-u_2^m|(x,t)&\leq C_T r^{-am}(\log(1/r))^{mq}
             &&\forall 0<r=|x|<r_T,\label{eq:smallbounds}\\
|u_1-u_2|(x,t)&\leq C_T r^{-b},
 & |u_1^m-u_2^m|(x,t)&\leq C_T r^{-(n-2)}
             &&\forall r=|x|>R_T.\label{eq:largebounds}
\end{align}
Fix \(s\geq0\), \(0<\theta<\mu\), and write
\(\Psi=\Psi_{s,\theta}\). Choose $0<2\rho<r_T$ and
$R>R_T$. Let $0\le\chi_{\rho,R}\le1$ be a smooth radial cutoff
that vanishes for $r\le\rho$ or $r\ge2R$ and equals one for
$2\rho\le r\le R$. It can be chosen so that
\[
\begin{array}{lll}
 |\nabla\chi_{\rho,R}|\le C/\rho,&
 |\Delta\chi_{\rho,R}|\le C/\rho^2,& \rho<r<2\rho,\\
 |\nabla\chi_{\rho,R}|\le C/R,&
 |\Delta\chi_{\rho,R}|\le C/R^2,& R<r<2R.
\end{array}
\]
Since
\[
\begin{split}
 &\bigl\||u_1(\cdot,h)-u_2(\cdot,h)|-|u_{0,1}-u_{0,2}|\bigr\|_{L^1(B_{2R}\setminus B_{\rho})}
 \\
 &\hspace{1cm}\le
 \|u_1(\cdot,h)-u_{0,1}\|_{L^1(B_{2R}\setminus B_{\rho})}
 +\|u_2(\cdot,h)-u_{0,2}\|_{L^1(B_{2R}\setminus B_{\rho})}\longrightarrow0\quad\forall 2R>\rho>0\quad\mbox{ as }h\to 0,
\end{split}
\]
 multiplying \eqref{eq:kato} by $\chi_{\rho,R}\Psi$ and integrating over $(\mathbb{R}^n\setminus\{0\})\times (h,t)$, $0<h<t\le T$, and then letting $h\to 0$, by Lemma  \ref{lem:superharmonic} we get
\begin{align}\label{eq:cutoffinequality}
 \int_{\mathbb{R}^n\setminus\{0\}} |u_1-u_2|(x,t)\chi_{\rho,R}\Psi\dd x
 &\leq \int_{\mathbb{R}^n\setminus\{0\}} |u_{0,1}-u_{0,2}|(x)\chi_{\rho,R}\Psi\dd x \notag\\
 &\quad+\frac{n-1}{m}\int_0^t\int_{\mathbb{R}^n\setminus\{0\}} |u_1^m-u_2^m|
       \bigl(\Psi\Delta\chi_{\rho,R}
               +2\nabla\Psi\cdot\nabla\chi_{\rho,R}\bigr)
         \dd x\dd\tau\notag\\
 &\leq \int_{\mathbb{R}^n\setminus\{0\}} |u_{0,1}-u_{0,2}|(x)\chi_{\rho,R}\Psi\dd x \notag\\
 &\quad+\frac{n-1}{m}\int_0^t\int_{\rho<|x|<2\rho} |u_1^m-u_2^m|
       \bigl(\Psi\Delta\chi_{\rho,R}
               +2\nabla\Psi\cdot\nabla\chi_{\rho,R}\bigr)
         \dd x\dd\tau\notag\\
         &\quad+\frac{n-1}{m}\int_0^t\int_{R<|x|<2R} |u_1^m-u_2^m|
       \bigl(\Psi\Delta\chi_{\rho,R}
               +2\nabla\Psi\cdot\nabla\chi_{\rho,R}\bigr)
         \dd x\dd\tau\notag\\
         &\leq \int_{\mathbb{R}^n\setminus\{0\}} |u_{0,1}-u_{0,2}|(x)\chi_{\rho,R}\Psi\dd x +I_1+I_2\quad\forall 0<t\le T. 
\end{align}
Here $I_1$ and $I_2$ denote the displayed space-time boundary integrals
over $\rho<|x|<2\rho$ and $R<|x|<2R$, including the factor $(n-1)/m$.
For \(\rho\le r\le 2\rho<1\), \(\Psi\leq C r^{-\mu+\theta}\).
By Lemma \ref{lem:superharmonic} and \eqref{eq:smallbounds},
\begin{align}
|I_1|\le&Ct\,\rho^{-2}\int_{\rho<|x|<2\rho}r^{-am}(\log(1/r))^{mq}r^{-\mu+\theta}\dd x\quad\forall 0<2\rho<r_T\notag\\
\le& C t\,\rho^{n-2-am-\mu+\theta}(\log(1/\rho))^{mq}\quad\forall 0<2\rho<r_T\notag\\
=&C t\,\rho^\theta(\log(1/\rho))^{m/(1-m)}
 \longrightarrow0\quad\mbox{ as }\rho\to 0. \label{eq:innererror}
\end{align}
The last equality uses \(am=a-2\) and \(\mu=n-a\).
If \(s>0\), then \(\Psi\leq C_s r^{-\mu-\delta+\theta}\) for large \(r\).
Since \(\theta<\mu\), hence 
\begin{equation}\label{eq:outererror}
|I_2|\le C_{T,s}t\,R^{-2}\!
 \int_{R<|x|<2R}r^{-(n-2)}r^{-\mu-\delta+\theta}\dd x
 \leq C_{T,s}t\,R^{-(\mu-\theta)-\delta}\longrightarrow0\quad\mbox{ as }R\to\infty.
\end{equation}
For \(s=0\), the same calculation gives the inner estimate \eqref{eq:innererror} and the outer estimate is
\begin{equation}\label{eq:outererror2}
 |I_2|\le C_T t R^{-\mu+\theta}\longrightarrow0\quad\mbox{ as }R\to\infty.
\end{equation}
For fixed $s\ge0$ and $\theta>0$, both cutoff integrands are
dominated by integrable functions. After including the radial volume
element, the bounds near zero and infinity are constant multiples of
\[
 r^{\theta-1}(\log(1/r))^q,\qquad
 r^{n-1-b-\mu+\theta},
\]
respectively. They are integrable because $\theta>0$, $b>n$, and
$\theta<\mu$; the same profile bounds apply to the initial data.
By the Lebesgue dominated convergence the cutoff bounds vanish as $\rho\to 0$ and $R\to\infty$.
Thus letting $\rho\to 0$ and $R\to\infty$ in 
\eqref{eq:cutoffinequality}, by \eqref{eq:innererror}, \eqref{eq:outererror} and \eqref{eq:outererror2} we get
\begin{equation}\label{eq:subcriticalcontraction}
 \int_{\mathbb{R}^n\setminus\{0\}} |u_1-u_2|(x,t)\Psi_{s,\theta}\dd x
 \leq\int_{\mathbb{R}^n\setminus\{0\}} |u_{0,1}-u_{0,2}|(x)\Psi_{s,\theta}\dd x.
\end{equation}
Let $0<\theta_0<\mu$. Then for any \(0<\theta\leq\theta_0\),
\begin{equation}\label{eq:difference-l1}
|u_{0,1}-u_{0,2}|(x)\Psi_{s,\theta}(x)\le G(x):=\left\{\begin{aligned}
&|u_{0,1}-u_{0,2}|(x)W_s(|x|)\quad\forall |x|\le 1\\
&|u_{0,1}-u_{0,2}|(x)|x|^{\theta_0}W_s(|x|)\quad\forall |x|>1.
\end{aligned}\right.
\end{equation}
By \eqref{eq:initial-data-weighted-difference}, 
\begin{equation}\label{eq:integral-finite1}
\int_{0<|x|\le 1} |u_{0,1}-u_{0,2}|(x)W_s\dd x<\infty.
\end{equation}
Since \(b>n\) and \(\theta_0<\mu\), by \eqref{eq:infinityasymp}, \eqref{eq:initial-value-upper-bd}  and \eqref{eq:difference-l1},
\begin{equation}\label{eq:integral-finite2}
\int_{|x|>1} |u_{0,1}-u_{0,2}|(x)|x|^{\theta_0}W_s(|x|)\dd x\le C\int_1^{\infty}r^{-1-(b-n)-(\mu-\theta_0)}\,dr<\infty.
\end{equation}
By \eqref{eq:difference-l1}, \eqref{eq:integral-finite1} and \eqref{eq:integral-finite2}, $G\in L^1(\mathbb{R}^n\setminus\{0\})$. Hence letting \(\theta\downarrow0\) through \(0<\theta\leq\theta_0<\mu\) in \eqref{eq:subcriticalcontraction}, by Fatou's lemma and Lebesgue dominated convergence theorem, we get
\eqref{eq:contraction}.

Similarly if \eqref{eq:initial-data-order} holds, then by a similar argument we get
\eqref{eq:contraction-order} and the proposition follows.
\end{proof}

We are now ready to prove Theorem \ref{thm:main}. 

\begin{proof}[Proof of Theorem \ref{thm:main}]
By Lemma \ref{lem:parameter-unique-and-range} the selected parameter $\lambda_0$ is unique and \(\lambda_2\le\lambda_0\le\lambda_1\).
Let \(S(t)=e^{-\beta t}\) and 
\[
 E_1(t)=\int_{\mathbb{R}^n\setminus\{0\}} |v(x,t)-f_{\lambda_0}(x)|W_1(|x|)\dd x.
\]
Since
\begin{equation}\label{eq:criticalidentity}
 \alpha-n\beta+\mu\beta=\beta(a-n+\mu)=0,
\end{equation}
by Proposition \ref{prop:contraction}, applied for the fixed value \(s=S(t)\) for
each fixed terminal time $t$. The parameter $s$ is held constant over
the entire physical-time interval $[0,t]$ in that application. We obtain
\begin{align}
 E_1(t)=&\int_{\mathbb{R}^n\setminus\{0\}}e^{\alpha t}
   |u(e^{\beta t}x,t)-U_{\lambda_0}(e^{\beta t}x,t)||x|^{-\mu}(1+|x|)^{-\delta}\dd x\notag\\
 =&e^{(\alpha-n\beta+\mu\beta) t}\int_{\mathbb{R}^n\setminus\{0\}}
   |u(y,t)-U_{\lambda_0}(y,t)||y|^{-\mu}(1+e^{-\beta t}|y|)^{-\delta}\dd y\notag\\
    =&\int_{\mathbb{R}^n\setminus\{0\}}
   |u(y,t)-U_{\lambda_0}(y,t)|W_{S(t)}(|y|)\dd y\notag\\
 \le&\int_{\mathbb{R}^n\setminus\{0\}} |u_0(y)-f_{\lambda_0}(y)|W_{S(t)}(|y|)\dd y. \label{eq:decayestimate}
\end{align}
Since \(\beta<0\), \(S(t)\to\infty\) as $t\to\infty$ and \(S(t)\geq 1\) for any $t>0$. For every \(r>0\),
\begin{equation}\label{eq:w-order}
 0<W_{S(t)}(r)\leq W_1(r)\leq r^{-\mu},
 \quad\mbox{ and }\quad W_{S(t)}(r)\to 0\quad\mbox{ as }t\to\infty.
\end{equation}
By \eqref{eq:mainhyp}, \eqref{eq:w-order} and the Lebesgue dominated convergence theorem the right hand side of 
\eqref{eq:decayestimate} tends to zero as $t\to\infty$. Thus by \eqref{eq:decayestimate},
\begin{equation}\label{eq:E1zero}
 E_1(t)\longrightarrow0\quad\mbox{ as }t\to\infty.
\end{equation}
By \eqref{eq:rescaed_function_trap}, \eqref{eq:lambda-ordering} and \eqref{eq:ordering},
\[
 |v(x,t)-f_{\lambda_0}(x)|
          \leq f_{\lambda_2}(x)-f_{\lambda_1}(x)\quad\forall x\in\mathbb{R}^n\setminus\{0\}.
\]
By \eqref{eq:infinityasymp}, for all sufficiently large \(R\),
\begin{align}
 \int_{|x|>R}|v-f_{\lambda_0}|\,|x|^{-\mu}\dd x
 &\leq C\int_R^\infty r^{n-1-b-\mu}\dd r\notag\\
 &\leq C R^{-(b-a)}
   =C R^{-\gamma_1}, \label{eq:uniformtail}
\end{align}
uniformly in \(t\). On \(0<r\leq R\),
\(r^{-\mu}\leq(1+R)^\delta W_1(r)\). Hence
\begin{equation}\label{eq:criticalupgrade}
 \int_{\mathbb{R}^n\setminus\{0\}} |v-f_{\lambda_0}|\,|x|^{-\mu}\dd x
 \leq (1+R)^\delta E_1(t)+C R^{-\gamma_1}.
\end{equation}
Letting first \(t\to\infty\) and then \(R\to\infty\) in \eqref{eq:criticalupgrade}, by \eqref{eq:E1zero} we get
 \eqref{eq:mainL1}. 

We will now prove the local \(C^2\) convergence.
We choose open annuli $A,A',A''$
such that $K\subset A$, $\overline A\subset A'$, $\overline{A'}\subset A''$,
and $\overline{A''}\Subset\mathbb R^n\setminus\{0\}$.
Let
\(z=v^m\). Then by \eqref{eq:rescaledpde},
\begin{equation}\label{eq:zPDE}
 z_t=(n-1)z^{(m-1)/m}\Delta z
          +\beta x\cdot\nabla z+m\alpha z\quad\mbox{ in }(\mathbb{R}^n\setminus\{0\})\times (0,\infty).
\end{equation}
Then by \eqref{eq:rescaed_function_trap}, there exists constants $C_{A''}>c_{A''}>0$ such that
\[
c_{A''}\leq v(x,t)\leq C_{A''}
 \quad\forall x\in A'',\ t\geq0.
\]
Then the equation \eqref{eq:zPDE} for $z$ is uniformly parabolic on $\overline{A''}$. The drift $\beta x$ and zeroth-order coefficient $m\alpha$ are also
bounded there. Apply the interior parabolic H\"older estimate for bounded measurable
coefficients~\cite{KS} on cylinders of fixed time length contained in
$A''\times(0,\infty)$. It gives a uniform parabolic H\"older bound
for $z$ on smaller cylinders in
$\overline{A'}\times[1/2,\infty)$.
Since $z$ remains in a fixed positive compact interval, the smooth
map $z\mapsto(n-1)z^{(m-1)/m}$ transfers that bound to the principal
coefficient. Interior Schauder estimates~\cite{LSU}, followed by
the smooth composition $v=z^{1/m}$, then give
\begin{equation*}
 \sup_{t\ge1}\|v(\cdot,t)\|_{C^{2+\eta}(\overline A)}<\infty
\end{equation*}
for some constant $\eta\in(0,1)$. Fixed cylinder lengths ensure that the
constants are independent of time translation.
 Then for any sequence $\{t_j\}\subset\mathbb{R}^+$, $t_j\to\infty$ as $j\to\infty$, the sequence $\{v(\cdot,t_j)\}_{j=1}^{\infty}$ is uniformly bounded and equi-H\"older continuous in $C^2(\overline A)$. Hence it has a subsequence which we may assume without loss of generality to be the sequence itself that converges in $C^2(A)$ as $j\to\infty$.
By  \eqref{eq:mainL1}  this subsequential limit is
\(f_{\lambda_0}\). Since the sequence is arbitrary, we get
\eqref{eq:mainC2}.
\end{proof}

\section{Convergence result for solutions without the lower profile bound}
\setcounter{equation}{0}
\setcounter{theorem}{0}

In this section we will prove Theorem \ref{thm:upperonly}. 

\begin{lemma}\label{lem:harnack-ineqn}
Let $n\ge 3$, $A\ge 0$ and \(D\subset\R^n\) be a connected open set. 
Then there exist constants  $C_0=C_0(n)>0$ and $r_0=r_0(A)>0$ such that 
for any function \(0\le z\in C^2(D)\) satisfying
\begin{equation}\label{eq:hypotheses}
\Delta z\leq Az\quad\text{in }D
\end{equation}
we have
\begin{equation}\label{eq:harnack}
 \inf_{B_r(x_0)} z\geq \frac{1}{C_0r^n}\int_{B_{2r}(x_0)} z(y)\dd y\quad\forall \overline{B_{5r}(x_0)}\subset D, 0<r\le r_0.
\end{equation}
\end{lemma}
\begin{proof}
The lemma follows immediately by applying Theorem 8.18 of \cite{GT} to
$\Delta-A$ with $p=1$. Alternatively the lemma also follows from Theorem 5.2 equation (5.6) of \cite{T} and H\"older's inequality.
For the sake of completeness we will give a simple elementary proof of \eqref{eq:harnack} here.
We divide the proof into two steps.

\medskip
\noindent\textbf{Step 1: an estimate at the center of one ball.}
Let \(x\in D\) and choose \(R>0\) such that
\begin{equation}\label{eq:Rchoice}
 \overline{B_R(x)}\subset D,\qquad \sqrt A\,R\leq1.
\end{equation}
Let
\begin{equation}\label{eq:cosineLift}
 Z(y,\xi):=z(y)\cos(\sqrt A\,\xi)\quad\forall y\in D, \xi\in\mathbb{R}.
\end{equation}
Then by \eqref{eq:Rchoice},
\begin{align*}
&|\sqrt A\,\xi|\leq1<\frac{\pi}{2}\quad\forall (x,\xi)\in B_R^{n+1}((x,0))\\
\Rightarrow\quad&Z(x,\xi)\geq0\quad\forall (x,\xi)\in B_R^{n+1}((x,0)).
\end{align*}
By \eqref{eq:hypotheses} and differentiation of \eqref{eq:cosineLift},
\begin{align*}
 \Delta_{y,\xi}Z
 &=\cos(\sqrt A\,\xi)(\Delta z-Az)\leq0\quad\forall (x,\xi)\in B_R^{n+1}((x,0)).
\end{align*}
Hence $Z$ is superharmonic on $B_R^{n+1}((x,0))$.
Thus by the  superharmonic mean-value inequality,
\begin{equation}\label{eq:meanvalue}
 z(x)=Z(x,0)
 \geq\frac{(n+1)}{\omega_{n+1}R^{n+1}}
 \int_{B_R^{n+1}((x,0))}Z(y,\xi)\dd y\dd\xi.
\end{equation}
Note that since
\[
 |y-x|^2+\xi^2<\frac{R^2}{4}+\frac{R^2}{4}<R^2\quad\forall (y,\xi)\in B_{R/2}^{n}(x)\times(-R/2,R/2),
\]
we have
\begin{equation}\label{eq:set-inclusion}
B_{R/2}^{n}(x)\times(-R/2,R/2)\subset B_R^{n+1}((x,0)).
\end{equation}
Moreover by \eqref{eq:Rchoice},
\begin{equation}\label{eq:cos-lower-bd}
\cos(\sqrt A\,\xi)\geq\cos(1/2)>0\quad\forall (y,\xi)\in B_{R/2}^{n}(x)\times(-R/2,R/2),
\end{equation}
Hence by \eqref{eq:meanvalue}, \eqref{eq:set-inclusion} and \eqref{eq:cos-lower-bd},
\begin{align*}
 z(x)
 &\geq\frac{(n+1)\cos(1/2)}{\omega_{n+1}R^{n+1}}
       \int_{B_{R/2}(x)}\int_{-R/2}^{R/2}z(y)\dd\xi\dd y\\
 &=\frac{(n+1)\cos(1/2)}{\omega_{n+1}R^n}
       \int_{B_{R/2}(x)}z(y)\dd y
\end{align*}
Thus
\begin{equation}\label{eq:pointmass}
z(x)\geq c_nR^{-n}\int_{B_{R/2}(x)}z(y)\dd y,
 \qquad c_n:=(n+1)\frac{\cos(1/2)}{\omega_{n+1}}>0.
\end{equation}

\noindent\textbf{Step 2: the integral over the doubled ball.}
Let
\[
 r_0=\frac{1}{4\sqrt{1+A}}.
\]
Then $r_0>0$. Suppose $0<r\le r_0$ and
$\overline{B_{5r}(p)}\subset D$. Then
\begin{equation}\label{eq:rchoiceLocal}
 4r\sqrt A\le 1,\qquad
 \overline{B_{4r}(x)}\subset B_{5r}(p)\quad\forall x\in B_r(p).
\end{equation}
For any $x\in B_r(p)$, 
\[
 \mathcal C=B_{2r}^n(p)\times(-r,r)\subset B_{4r}^{n+1}((x,0)).
\]
Indeed, for $(y,\xi)\in\mathcal C$,
\[
 |y-x|<3r,\quad |\xi|<r,
 \qquad |y-x|^2+\xi^2<10r^2<16r^2.
\]
On this cylinder $|\sqrt A\,\xi|\le1/4$, so
$\cos(\sqrt A\,\xi)\ge\cos(1/2)>0$.
Apply \eqref{eq:meanvalue} with radius $4r$ and restrict its integral
to $\mathcal C$. The $\xi$ interval has length $2r$, and hence
\[
 z(x)\ge
 \frac{2(n+1)\cos(1/2)}{4^{n+1}\omega_{n+1}}\,
 r^{-n}\int_{B_{2r}(p)}z(y)\dd y.
\]
Taking the infimum over $x\in B_r(p)$ proves
\begin{equation*}
 \inf_{B_r(p)}z\ge b_n r^{-n}\int_{B_{2r}(p)}z(y)\dd y,
 \qquad b_n=\frac{2(n+1)\cos(1/2)}{4^{n+1}\omega_{n+1}}>0.
\end{equation*}
Thus \eqref{eq:harnack} holds with $C_0=b_n^{-1}$.

\end{proof}

\begin{lemma}[A lower bound from positive mass]\label{lem:masspositivity}
Let $n\ge 3$, \(D\subset\R^n\) be a connected open set, \(K\subset D\) be a nonempty compact set, and \(E\subset D\) be a compact measurable set of positive measure.
Let \(A\geq0\) be a constant. Then there exists a constant \(c=c(n,D,K,E,A)>0\) such that
for any function \(0\le z\in C^2(D)\) satisfying
\eqref{eq:hypotheses} 
we have
\begin{equation}\label{eq:masspositivity}
 \inf_K z\geq c\int_E z(y)\dd y.
\end{equation}
\end{lemma}
\begin{proof}
We divide the proof into three steps.

\medskip
\noindent\textbf{Step 1: construction of  a finite connected family of balls.}
For any $x\in E\cup K$, there exists an open ball $B_{r_x}(x)$ such that 
$\overline{B_{r_x}(x)}\subset D$. Since $E\cup K$ is compact, the family of open cover $\{B_{r_x}(x)\}_{x\in E\cup K}$ has a finite subcover
$\{B_{r_i}(x_i)\}_{i=1}^{n_0}$ such that $E\cup K\subset\cup_{i=1}^{n_0}B_{r_i}(x_i)$.
Since $D$ is an open connected set, $D$ is path connected. Hence for any $i=2,\cdots,n_0$,
there exists continuous path $\gamma_i:[0,1]\to D$ such that $\gamma_i(0)=x_1$, $\gamma_i(1)=x_i$. Let
\begin{equation*}
K_1=\left(\cup_{i=1}^{n_0}\overline{B_{r_i}(x_i)}\right)\cup\left(\cup_{i=2}^{n_0}\gamma_i([0,1])\right).
\end{equation*}
Then $K_1$ is a compact path connected set and $E\cup K\subset K_1\subset D$. Let $r_0>0$ and $C_0>0$ be given by Lemma \ref{lem:harnack-ineqn}. If \(D\neq\mathbb{R}^n\), then since $K_1$ is compact, \(\delta:=\dist(K_1,\R^n\setminus D)>0\). If \(D=\mathbb{R}^n\), let $\delta=\infty$. Let $\delta_0=\min (\delta, r_0)$. Then for any 
\(0<r<\delta_0/5\),
\begin{equation*}
 \overline{B_{5r}(p)}\subset D\qquad\text{for every }p\in K_1.
\end{equation*}
Since \(K_1\) is compact, we can cover \(K_1\) by finitely many balls with centers in \(K_1\):
\[
K_1\subset\bigcup_{i=1}^N U_i,
 \qquad U_i:=B_{r/2}(p_i),\qquad p_i\in K_1.
\]
Set \(B_i:=B_r(p_i)\), then \(\{B_i\}_{i=1}^N\) cover \(E\cup K\).

Make a graph whose vertices are \(p_1,\ldots,p_N\), joining \(p_i\) and
\(p_j\) when \(K_1\cap U_i\cap U_j\neq\varnothing\). This graph is
connected: otherwise two nonempty groups of graph components would
give disjoint nonempty relatively open subsets covering $K_1$, since
every ball meets $K_1$ at its center. That would contradict the
connectedness of $K_1$.
If \(p_i,p_j\) are adjacent, then \(K_1\cap U_i\cap U_j\neq\varnothing\). Let \(p\in K_1\cap U_i\cap U_j\). Then
\begin{align}\label{eq:centers}
 &|p_i-p_j|\le |p_i-p|+|p-p_j|<r\notag\\
\Rightarrow\quad&|y-p_j|\leq|y-p_i|+|p_i-p_j|<2r\quad\forall y\in B_i\notag\\
\Rightarrow\quad& B_i\subset B_{2r}(p_j).
\end{align}
Note that any two vertices of a connected graph with \(N\) vertices can be
joined by a simple path containing at most \(N-1\) edges.

\medskip
\noindent\textbf{Step 2: an initial lower bound.}
Let \(0\le z\in C^2(D)\) satisfy \eqref{eq:hypotheses}. Let
\[
 m_i:=\inf_{B_i}z,
 \qquad M:=\int_E z(y)\dd y.
\]
If \(M=0\), \eqref{eq:masspositivity} holds for any constant $c>0$.
Suppose \(M>0\). Since \(E\subset\bigcup_{i=1}^N B_i\) and \(z\geq0\),
\[
 M\leq\sum_{i=1}^N\int_{B_i}z(y)\dd y.
\]
Thus there is an index \(i_0\) for which
\(\int_{B_{i_0}}z\geq M/N\). By Lemma \ref{lem:harnack-ineqn},
\begin{equation}\label{eq:HarnackSeed}
 m_{i_0}
 \geq\frac{1}{C_0r^n}\int_{B_{2r}(p_{i_0})}z(y)\dd y
 \geq\frac{M}{C_0Nr^n}.
\end{equation}

\medskip
\noindent\textbf{Step 3: propagate the bound along the graph.}
For adjacent vertices \(i,j\), by \eqref{eq:centers} and
\eqref{eq:harnack} at \(p_j\),
\begin{align}
 m_j
 &\geq\frac{1}{C_0r^n}\int_{B_{2r}(p_j)}z(y)\dd y\notag\\
 &\geq\frac{1}{C_0r^n}\int_{B_i}z(y)\dd y
 \geq\frac{|B_i|}{C_0r^n}m_i
 =\frac{\omega_n}{nC_0}m_i.
 \label{eq:HarnackPropagation}
\end{align}
Set
\[
 \theta:=\min\{1,\omega_n/(nC_0)\}>0.
\]
For any \(i\), choose a simple path from \(i_0\) to \(i\), of length
\(\ell\leq N-1\). By iteration of \eqref{eq:HarnackPropagation} and
\eqref{eq:HarnackSeed},
\[
 m_i\geq\theta^\ell m_{i_0}
 \geq\frac{\theta^{N-1}}{C_0Nr^n}M.
\]
Since  $\{B_i\}_{i=1}^N$
cover \(K\), we conclude that
\begin{equation*}
\inf_K z\geq
 \frac{\theta^{N-1}}{C_0Nr^n}\int_E z(y)\dd y
\end{equation*}
and \eqref{eq:masspositivity} follows.

\end{proof}

We are now ready for the proof of Theorem \ref{thm:upperonly}.

\begin{proof}[Proof of Theorem \ref{thm:upperonly}]
\textbf{Step 1: selection of the parameter and approximation of the data.}
By the same argument as the proof of Lemma \ref{lem:parameter-unique-and-range}, $\lambda_0$ is unique.
If \(\lambda_0<\bar\lambda\), then by \eqref{eq:upperinitial},
\[
 f_{\lambda_0}-u_0
 \geq f_{\lambda_0}-f_{\bar\lambda}>0\quad\mbox{ on }\mathbb{R}^n\setminus\{0\}.
\]
This together with Lemma \ref{lem:profiledistance} implies that the left hand side integral of \eqref{eq:mainhyp} is infinite which contradicts \eqref{eq:mainhyp}. Hence
\(\lambda_0\geq\bar\lambda\). 

We choose an increasing sequence $\{L_j\}_{j=1}^{\infty}\subset (\lambda_0,\infty)$, \(L_j\to\infty\) as $j\to\infty$, and let
\begin{equation*}
 u_{0,j}(x)=\max\{u_0(x),f_{L_j}(x)\}.
\end{equation*}
Then by \eqref{eq:ordering},
\begin{equation}\label{u0j-compare2}
 f_{L_j}\leq u_{0,j}\leq f_{\bar\lambda},\qquad
 u_{0,j+1}\leq u_{0,j}\quad\mbox{ in }\mathbb{R}^n\setminus\{0\}\quad\forall j\in\mathbb{Z}^+.
\end{equation}
For each fixed $x\ne0$, the scaling formula and the far-field
asymptotic give
\[
 f_{L_j}(x)=L_j^a f_1(L_jx)
 \sim L_j^{a-b}|x|^{-b}\longrightarrow0,
\]
because $b>a$. Therefore
\begin{equation}\label{u0j-monone-decreasing}
u_{0,j}\downarrow u_0\quad\mbox{  a.e. in $\mathbb{R}^n\setminus\{0\}$ as }j\to\infty.
\end{equation}
Since 
\begin{equation*}
f_{L_j}\leq f_{\lambda_0}\quad\mbox{ in }\mathbb{R}^n\setminus\{0\}\quad\forall j\in\mathbb{Z}^+,
\end{equation*}
 we have
\begin{align}
 0\leq u_{0,j}-u_0
 &\leq(f_{\lambda_0}-u_0)_+
        \leq |u_0-f_{\lambda_0}|,\quad\mbox{ in }\mathbb{R}^n\setminus\{0\}\quad\forall j\in\mathbb{Z}^+\label{eq:upperapproxerror}\\
 |u_{0,j}-f_{\lambda_0}|
 &\leq |u_0-f_{\lambda_0}|.\quad\mbox{ in }\mathbb{R}^n\setminus\{0\}\quad\forall j\in\mathbb{Z}^+.\label{eq:upperapproxdistance}
\end{align}
Then by \eqref{eq:mainhyp}, \eqref{u0j-monone-decreasing}, \eqref{eq:upperapproxerror} and the Lebesgue dominated convergence theorem,
\begin{equation}\label{eq:dj}
 d_j:=\int_{\mathbb{R}^n\setminus\{0\}}(u_{0,j}(x)-u_0(x))|x|^{-\mu}\dd x\longrightarrow0\quad\mbox{ as }j\to\infty.
\end{equation}

\textbf{Step 2: the trapped approximating solutions and their estimates.}
By Theorem 1.4 of \cite{HP} for any \(j\in\mathbb{Z}^+\) there exists  a positive classical
solution \(u_j\) of \eqref{eq:pde} with initial datum \(u_{0,j}\) satisfying
\begin{equation}\label{eq:ujtrapAB}
 U_{L_j}\leq u_j\leq U_{\bar\lambda},
 \qquad
 (u_j)_t\leq\frac{u_j}{(1-m)t}\quad\mbox{ in }(\mathbb{R}^n\setminus\{0\})\times (0,\infty).
\end{equation}
Then by \eqref{u0j-compare2}, \eqref{eq:ujtrapAB} and Proposition \ref{prop:contraction},
\begin{equation}\label{eq:uj-comparison-ineqn}
u_{j+1}\leq u_j\quad\mbox{ in }(\mathbb{R}^n\setminus\{0\})\times (0,\infty)\quad\forall j\in\mathbb{Z}^+.
\end{equation}
Let
\[
 v_j(x,t)=e^{\alpha t}u_j(e^{\beta t}x,t)\quad\forall x\in \mathbb{R}^n\setminus\{0\},t>0.
\]
Then by \eqref{eq:ujtrapAB} and \eqref{eq:uj-comparison-ineqn},
\begin{equation}\label{eq:vjupper}
 0<v_{j+1}\leq v_j\leq f_{\bar\lambda}\quad\mbox{ in }(\mathbb{R}^n\setminus\{0\})\times (0,\infty)\quad\forall j\in\mathbb{Z}^+.
\end{equation}
By  Proposition \ref{prop:contraction}, \eqref{eq:criticalidentity} and
\eqref{eq:upperapproxdistance},
\begin{equation}\label{eq:vjcritical}
 \int_{\mathbb{R}^n\setminus\{0\}} |v_j(x,t)-f_{\lambda_0}(x)||x|^{-\mu}\dd x
 \leq \int_{\mathbb{R}^n\setminus\{0\}} |u_{0,j}(x)-f_{\lambda_0}(x)||x|^{-\mu}\dd x
 \leq E_0
 \quad\forall j\geq1,\ t>0.
\end{equation}
By Theorem \ref{thm:main} for any $j\in\mathbb{Z}^+$,
\begin{equation}\label{eq:mainL1-j}
\lim_{t\to\infty}\int_{\mathbb{R}^n\setminus\{0\}}
       |v_j(x,t)-f_{\lambda_0}(x)|\,|x|^{-\mu}\dd x=0
\end{equation}
and
\begin{equation}\label{eq:mainC2-j}
 \lim_{t\to\infty}\norm{v_j(\cdot,t)-f_{\lambda_0}}_{C^2(K)}=0
 \quad\hbox{for any compact set } K\subset\mathbb{R}^n\setminus\{0\}. 
\end{equation}

\textbf{Step 3: a positive amount of mass is present on one fixed annulus.}
The leading profile asymptotic \eqref{eq:originasymp} for $f_{\lambda_0}$ and \(\mu+a=n\) imply
\begin{equation*}
 \int_{0<|x|<r_0} f_{\lambda_0}(x)|x|^{-\mu}\dd x=\infty,
\end{equation*}
for any $r_0>0$ since its radial integrand is comparable to
\(r^{-1}(\log(1/r))^{1/(1-m)}\).
Then we can choose a closed annulus
\[
 A_*=\{r_*^-\leq |x|\leq r_*^+\}\subset{\mathbb{R}^n\setminus\{0\}},\quad r_*^+>r_*^->0
\]
with $r_*^-$ sufficiently small such that
\begin{equation}\label{eq:anchorchoice}
 \int_{A_*}f_{\lambda_0}(x)|x|^{-\mu}\dd x>E_0+1.
\end{equation}
By \eqref{eq:vjcritical} and \eqref{eq:anchorchoice},
\begin{equation}\label{eq:anchormass}
 \int_{A_*}v_j(x,t)|x|^{-\mu}\dd x\geq1
 \quad\forall  j\geq1,\ t>0.
\end{equation}

\textbf{Step 4: convert the time-derivative bound into local positivity.}
Let \(K\subset{\mathbb{R}^n\setminus\{0\}}\) be a compact set and \(h>0\).
Choose a bounded connected smooth annular domain
\(D\subset{\mathbb{R}^n\setminus\{0\}}\) containing \(A_*\) and \(K\) in its interior.
Let $z_j=v_j^m$ for any $j\in\mathbb{Z}^+$ and 
\[
 M_D=\sup_{\overline D}f_{\bar\lambda}<\infty.
\]
Since \(0<v_j\leq M_D\) on \(D\) and \(m<1\),
\[
 z_j=v_j^m=v_jv_j^{m-1}\geq M_D^{m-1}v_j\quad\forall j\in\mathbb{Z}^+.
\]
Thus \eqref{eq:anchormass} implies
\begin{equation}\label{eq:anchorzmass}
 \int_{A_*}z_j(x,t)\dd x
 \geq (r_*^-)^{\mu}M_D^{m-1}>0\quad\forall j\in\mathbb{Z}^+,t>0.
\end{equation}
The identity \(m\alpha+2\beta=\alpha\) gives the exact spatial
scaling formula
\begin{align*}
 \Delta(v_j^m)(x,t)
 &=e^{(m\alpha+2\beta)t}
     \Delta(u_j^m)(e^{\beta t}x,t)\notag\\
 &=\frac{m}{n-1}e^{\alpha t}
                  (u_j)_t(e^{\beta t}x,t)\notag\\
 &\leq \frac{m}{(n-1)(1-m)t}\,v_j(x,t)\quad\mbox{ in }\mathbb{R}^n\setminus\{0\}\quad\forall j\in\mathbb{Z}^+,t>0.
\end{align*}
For \(t\geq h\), it implies
\begin{equation*}
 \Delta z_j\leq A_{D,h}z_j\quad\mbox{ in }D,\quad
 A_{D,h}=\frac{mM_D^{1-m}}{(n-1)(1-m)h},\quad\forall j\in\mathbb{Z}^+,t\ge h
\end{equation*}
because \(v_j=z_jv_j^{1-m}\leq M_D^{1-m}z_j\).
Apply Lemma \ref{lem:masspositivity} to \(z_j(\cdot,t)\) on \(D\),
with \(E=A_*\), by \eqref{eq:anchorzmass},
\begin{equation}\label{eq:vjuniformlower}
 \inf_{\substack{x\in K\\ t\geq h}}v_j(x,t)
 \geq c_{K,h}>0
 \quad\forall j\in\mathbb{Z}^+.
\end{equation}

\textbf{Step 5: construction of  a positive classical limit.}
Since \(K\) and \(h\) are arbitrary, by
\eqref{eq:vjupper} and \eqref{eq:vjuniformlower} the sequence \(\{v_j\}_{j=1}^{\infty}\) is uniformly bounded above and below away from zero on every compact space-time cylinder
contained in \((\mathbb{R}^n\setminus\{0\})\times(0,\infty)\).
Hence the equation \eqref{eq:zPDE} for the sequence \(\{z_j\}_{j=1}^{\infty}\),  is uniformly parabolic on those cylinders. By an argument similar to the proof of
Theorem~\ref{thm:main} we get uniform interior H\"older \cite{KS} and Schauder estimates \cite{LSU} on every compact subset of 
\((\mathbb{R}^n\setminus\{0\})\times(0,\infty)\). 
These then give local uniform
$C^{2+\eta,1+\eta/2}$ bounds, first for $z_j$ and then for $v_j$,
 on every compact subset of 
\((\mathbb{R}^n\setminus\{0\})\times(0,\infty)\).
The monotone sequence $v_j$ has a pointwise limit $v$,
which is positive by \eqref{eq:vjuniformlower}.
The local estimates and a diagonal compactness argument imply
\[
 v_j\longrightarrow v,\qquad z_j\longrightarrow v^m
 \quad\hbox{in }C^{2,1}_{\loc}
 \bigl((\mathbb R^n\setminus\{0\})\times(0,\infty)\bigr)\quad\mbox{ as }j\to\infty.
\]
Passing to the limit in the rescaled equation shows that $v$ satisfies
\eqref{eq:rescaledpde}.
Its inverse rescaling
\[
 u(x,t)=e^{-\alpha t}v(e^{-\beta t}x,t)
\]
is a positive classical solution of \eqref{eq:fde} in $(\mathbb{R}^n\setminus\{0\})\times(0,\infty))$, is the decreasing
limit of \(u_j\), and satisfies the upper bound in
\eqref{eq:upperclass}. The derivative inequality in
\eqref{eq:ujtrapAB} passes to the local classical limit and gives
\eqref{eq:upperAB}. Equation \eqref{eq:vjuniformlower} also passes
to the limit, proving \eqref{eq:upperlocalpositive}.

\textbf{Step 6: the initial trace and uniqueness.}
For any \(\ell>j\ge 1\), by Proposition \ref{prop:contraction} and monotonicity,
\begin{align*}
 \int_{\mathbb{R}^n\setminus\{0\}} |v_j(x,t)-v_\ell(x,t)||x|^{-\mu}\dd x
 &= \int_{\mathbb{R}^n\setminus\{0\}} |u_j(x,t)-u_\ell(x,t)||x|^{-\mu}\dd x\notag\\
 &\le\int_{\mathbb{R}^n\setminus\{0\}} (u_{0,j}(x)-u_{0,\ell}(x))|x|^{-\mu}\dd x\notag\\
 &\leq d_j. 
\end{align*}
Letting \(\ell\to\infty\) and using Fatou's lemma gives
\begin{equation}\label{eq:limitstability}
 \int_{\mathbb{R}^n\setminus\{0\}} |u_j(x,t)-u(x,t)||x|^{-\mu}\dd x
 =
 \int_{\mathbb{R}^n\setminus\{0\}} |v_j(x,t)-v(x,t)||x|^{-\mu}\dd x
 \leq d_j
 \quad\forall t>0.
\end{equation}
For any non-empty compact set \(K\subset{\mathbb{R}^n\setminus\{0\}}\), let
$R_K=\max_{x\in K}|x|$. Since $|x|^{-\mu}\ge R_K^{-\mu}$ on $K$,
by \eqref{eq:dj} and \eqref{eq:limitstability},
\begin{align}
 \norm{u(\cdot,t)-u_0}_{L^1(K)}
 &\leq
 \norm{u(\cdot,t)-u_j(\cdot,t)}_{L^1(K)}
 +\norm{u_j(\cdot,t)-u_{0,j}}_{L^1(K)}
 +\norm{u_{0,j}-u_0}_{L^1(K)}\notag\\
 &\leq 2R_K^{\mu}d_j
        +\norm{u_j(\cdot,t)-u_{0,j}}_{L^1(K)}.
        \label{eq:upperinitialtrace}
\end{align}
Letting first $t\to\infty$ and then $j\to\infty$ in \eqref{eq:upperinitialtrace} we get that
$u$ has initial trace $u_0$. If $\widehat u$ is another positive classical solution with initial trace $u_0$ which satisfies $0<\widehat u\le U_{\bar\lambda}$ in $(\mathbb{R}^n\setminus\{0\})\times (0,\infty)$, then by Proposition \ref{prop:contraction} $u=\widehat u$ in $(\mathbb{R}^n\setminus\{0\})\times (0,\infty)$. This proves uniqueness in the
upper-envelope class stated in the theorem.

\textbf{Step 7: weighted convergence at large times.}
By the triangle inequality and \eqref{eq:limitstability},
\[
 \norm{v(\cdot,t)-f_{\lambda_0}}_{L^1(|x|^{-\mu})}
 \leq d_j+\norm{v_j(\cdot,t)-f_{\lambda_0}}_{L^1(|x|^{-\mu})}\quad\forall j\ge 1.
\]
Letting first \(t\to\infty\)
and then \(j\to\infty\), by \eqref{eq:dj} and \eqref{eq:mainL1-j} we get
\eqref{eq:mainL1}.

\textbf{Step 8: local \(C^2\) convergence.}

By the bound \eqref{eq:upperlocalpositive} on nested annuli and for
$t\ge1/2$ and an argument similar to the proof of Theorem~\ref{thm:main} we get interior H\"older estimates \cite{KS},
followed by Schauder estimates \cite{LSU} and therefore
uniform $C^{2+\eta}$ bounds for $v(\cdot,t)$ on each smaller
closed annulus for $t\ge1$ where $\eta\in (0,1)$ is some constant.
Every sequence \(t_k\to\infty\) thus has a subsequence converging
in \(C^2\) on an annulus containing any prescribed compact \(K\).
The weighted convergence just proved identifies its limit there
with \(f_{\lambda_0}\). This proves \eqref{eq:mainC2} and completes the proof of the theorem.

\end{proof}


\begin{thebibliography}{99}


\bibitem[A]{A}
D.~G. Aronson, {\it The porous medium equation}, in
\emph{Nonlinear Diffusion Problems} (Montecatini Terme, 1985),
Lecture Notes in Mathematics, vol.~1224, Springer, Berlin, 1986, pp.~1--46.


\bibitem[DK]{DK}
P. Daskalopoulos and C.~E. Kenig,
{\it Degenerate Diffusions: Initial Value Problems and Local Regularity Theory},
EMS Tracts in Mathematics, vol.~1, European Mathematical Society,
Z\"urich, 2007. 


\bibitem[DS1]{DS1} P.~Daskalopoulos and N.~Sesum, {\em On the extinction profile of solutions to fast diffusion}, J. Reine Angew. Math. 622 (2008), 95--119.

\bibitem[DS2]{DS2} P.~Daskalopoulos and N.~Sesum, {\em The classification of locally conformally flat Yamabe solitons}, Advances in Math. 240 (2013), 346-369.

\bibitem[GT]{GT} D.~Gilbarg and N.~S. Trudinger, {\it Elliptic Partial Differential Equations of Second Order}, second edition, revised third printing (1998), Classics in Mathematics reprint, Springer, Berlin, 2001.

\bibitem[Hs]{Hs} S.Y. Hsu, {\it Existence and asymptotic behaviour of solutions of the very fast diffusion equation}, Manuscripta Math. 140 (2013), no. 3-4, 441-460.



\bibitem[H1]{H1} K.M.~Hui, {\it Asymptotic behaviour of solutions of the fast diffusion equation near its extinction time}, J. Math. Anal. Appl. 454 (2017), Issue 2,  695--715.

\bibitem[H2]{H2} K.M.~Hui, {\it Uniqueness and time oscillating behaviour of finite points blow-up solutions of the fast diffusion equation}, Proceedings of the Royal Society of Edinburgh Section A: Mathematics  150 (2020), no. 6, 2849--2870.




\bibitem[HP]{HP} K.M.~Hui and Jinwan Park, {\it Asymptotic behaviour of singular solution of the fast diffusion equation in the punctured euclidean space}, Discrete Contin. Dyn. Syst. 41 (2021), Issue 11, 5473-5508. 

\bibitem[K]{K}  T. Kato, {\it Schr\"odinger operators with singular potentials},
Israel J. Math., 13 (1972), 135--148.

\bibitem[KS]{KS}
N.~V. Krylov and M.~V. Safonov,
{\it A certain property of solutions of parabolic equations with
measurable coefficients}, Math. USSR-Izv. 16 (1981), no.~1,
151--164.

\bibitem[LSU]{LSU}
O.~A. Ladyzhenskaya, V.~A. Solonnikov and N.~N. Ural'tseva,
{\it Linear and Quasilinear Equations of Parabolic Type},
Translations of Mathematical Monographs, vol.~23,
American Mathematical Society, Providence, RI, 1968
(English translation of the Russian original).

\bibitem[T]{T} N.~S. Trudinger, {\it Linear elliptic operators with measurable coefficients}, Annali della Scuola Normale Superiore di Pisa, Classe di Scienze,
series~3, 27 (1973), no.~2, 265--308.


\bibitem[V]{V} J.L. Vazquez, {\it Smoothing and Decay Estimates for Nonlinear Diffusion Equations}, Oxford Lecture Series in Mathematics and its Applications 33, Oxford University Press, Oxford, 2006.




\end{thebibliography}
\end{document}